\documentclass{amsart}
\usepackage{amsmath,amsthm,amsfonts,latexsym,amssymb}
\usepackage{enumitem, blindtext}
\makeatletter
\def\@settitle{\begin{center}%
  \baselineskip14\p@\relax
  \bfseries
  \uppercasenonmath\@title
  \@title
  \ifx\@subtitle\@empty\else
     \\[1ex]\uppercasenonmath\@subtitle
     \footnotesize\mdseries\@subtitle
  \fi
  \end{center}%
}
\def\subtitle#1{\gdef\@subtitle{#1}}
\def\@subtitle{}
\makeatother
\usepackage{graphicx}
\newtheorem{theorem}{Theorem}

\newtheorem*{lemma}{Lemma}

\theoremstyle{definition}
\newtheorem{proposition}[theorem]{Proposition}

\newtheorem*{remark}{Remark}

\numberwithin{equation}{section}

\begin{document}

\title{On group rings}

%    Remove any unused author tags.

%    author one information
\author{Guy Renault}
\address{}
\thanks{Translator's note:
The original publication in French was ``Sur les anneaux de groupes.'' CR Acad. Sci. Paris Sér. AB 273 (1971): A84-A87.  This translation was prepared by Ryan C. Schwiebert using Google Translator Toolkit, and was graciously proofread by R. Christopher Coski, professor of French at Ohio University. Where mathematician and machine made mistakes, Chris Coski concocted corrections. Thanks also to the Acad\'emie des Sciences for authorizing the publication of this translation to arXiv.org in April 2018.}

\date{Session of June 28, 1971. This translation October 2017}

\dedicatory{Note by M. Guy Renault,\\ 
presented by M. Jean Leray.}

\begin{abstract}

We characterize the rings $A$ and groups $G$ for which the group rings
$A [G]$ are local, semi-local, or left perfect \cite{faculte1970seminaire}. The recent work of
M. P. Malliavin \cite{malliavin1971anneaux} and J. L. Pascaud permits the completion of results of \cite{faculte1970seminaire}
on self-injective group rings. 
\end{abstract}

\maketitle

\vskip.25cm

$A$ designates a ring with identity but which is not necessarily commutative, and $G$ is a
group. The fields involved are not necessarily commutative. 
For an exposition on group rings, consult J. Lambek \cite{lambek1966lectures} and P. Ribenboim \cite{ribenboim1969rings}.

\section{Local group rings}

We generalize a result of
T. Gulliksen-P. Ribenboim-T. M. Viswanathan \cite[p.~153]{gulliksen1970elementary} obtained for
the class of commutative group rings.

\begin{theorem}\label{thm:1} Let $A$ be a ring and $G$ a group $\neq e$ such that the group ring $A [G]$ is local. We then have the following properties:
\begin{enumerate}[label=(\alph*)]
\item $A$ is a local ring whose maximal left ideal will be denoted by $M$.
\item The field $K \neq A / M$ has characteristic $ p \neq  0$.
\item $G$ is a $p$-group. 
\end{enumerate}

If, additionally, $G$ is locally finite, these conditions are sufficient for $A [G]$ to be local.
\end{theorem}

The ring $A$ is isomorphic to a quotient ring of $A [G]$, hence $(a)$.
For the same reason $K [G]$ is a local ring. If $H$ is a subgroup
of $G$, then $K [H]$ is local. Indeed, let $R$ (resp. $R'$) be the radical of $K [G]$
(resp.  $K [H])$. It follows from a result of Connell \cite[p.~665]{connell1963group} that $ K [H] \cap R \subset R' $;
since $R$ is is the fundamental ideal of $K [G]$, $K [H] \cap R$ is the fundamental ideal of $K [H]$: this is a maximal left ideal which is equal to $R'$ and $K [H]$ is local. Let $ x \neq e$ be an element of $G$, $H_0$ the subgroup generated by $x$. $K [H_0]$ is a local ring and consequently the element $ e +x-x^2 $ is invertible. It is easy to see that this last condition implies
the finiteness of $ H_0$. Let $q$ be the order of $x$. If $q$ is invertible in $K$, the element 
$ e - q ^ {- i} \sum_{i = 0} ^ {q-1}$ $ x_i$ would be a nontrivial idempotent of $K [G]$ which is
not possible. We deduce immediately properties $(b)$ and $(c)$.

Let $A$ be a ring, $G$ a locally finite group satisfying the conditions of the theorem. Since $G$ is locally finite, $MA [G]$ is contained in the
radical of  $A[G]$ \cite[p.~665]{connell1963group} and it is sufficient to demonstrate that the ring $K [G]$
is local, which easily results from the following property that is well-known
when the field is commutative. Let $K$ be a (not necessarily commutative) field of characteristic $p \neq 0$, and $G$ a finite $p$-group.
Then $K [G]$ is a local ring whose radical is a nil ideal.

\begin{remark}Let $G$ be the infinite $p$-group generated by three elements
that is described in \cite{herstein1971noncommutative}, and let $k$ be the field of $p$ elements. $k [G]$ is a local ring
although $G$ is not locally finite.
\end{remark}

In what follows, $A$ and $G$ are commutative. The result of \cite[p.~153]{gulliksen1970elementary}
can also be generalized in the following way:

\begin{proposition} \label{prop:2} The following conditions are equivalent:
\begin{enumerate}[label=(\arabic*)]
\item $A[G]$ is a semi-local ring.
\item 
\begin{enumerate}
\item $A$ is a semi-local ring with radical $R$;
\item $G$ is finite or $G$ is infinite and in this case $A / R$ is a ring
of characteristic $p \neq 0$, $G = G_p  \times  G_0$ where $G_p$ is an infinite $p$-group,  and where $G_0$
is a finite group whose order is not divisible by $p$.
\end{enumerate}
\end{enumerate}
\end{proposition}

The proof of this theorem is not difficult. For the implication
(1)$\implies$ (2) consult \cite{burgess1969semi}.

\section{Left perfect group rings \cite{bass1960finitistic}}

Let's recall that
if $A$ is left perfect, the finitely generated sub-modules of any right $A$-module
satisfy the descending chain condition \cite{bjork1969rings}. The result that
follows was also obtained by Sheila Woods \cite{woods1971perfect} by completely different methods.

\begin{theorem}\label{thm:3}Let $A$ be a ring and $G$ be a group. The following are equivalent:
\begin{enumerate}[label=(\arabic*)]
\item $A [G]$ is left perfect.
\item
\begin{enumerate}\item $A$ is left perfect.
\item $G$ is finite.
\end{enumerate}
\end{enumerate}
\end{theorem}

(2) $\implies$ (1): For the finitely generated right ideals of $A [G]$, which are
finitely generated right $A$-modules, we verify the descending chain condition \cite{bjork1969rings}.

(1) $\implies$ (2): Let $R$ be the radical of $A$. The rings $A$, $(A / R) [G]$, which are
quotient rings of $A [G]$, are left perfect and it is sufficient to study
the case when $A$ is a simple ring with center $k$. $A [G]$ is a free $k [G]$ -module
, Lemma 12 of \cite[p.~150]{ribenboim1969rings} and the results of \cite{bjork1969rings} show that $k [G]$
is left perfect. Suppose $G$ is infinite: then $k [G]$ is not semiprimary
and it results in the following consequences: the characteristic of $k$ is $p> 0$
and there is a normal subgroup of $H_1$ of $G$ whose order is divisible
  by $p$ \cite[p.~162]{lambek1966lectures}, $G / H_1$ is infinite and $k [G / H_1]$ is left perfect. There
exists a normal sub-group $H_2$ of $G$ containing $H_1$ such that $p$ divides
the order of $H_2 / H_1$. Evidentially, then, there is an increasing sequence of
normal subgroups ($H_n$) of $G$ of order $p^{s(n)}q_n$, $p$ not dividing $q_n$, such
that $ s (n)> s (n - 1) $. The Sylow theorems permit the construction of
an infinite strictly increasing sequence of finite $p$-groups whose union
is an infinite $p$-group  $G_0$. $k [G]$ is a free $k [G_0]$ -module by Lemma 12
of \cite[p.~150]{ribenboim1969rings} and the results of \cite{bjork1969rings} show that $k [G_0]$ is left perfect.
$k [G_0]$ is a local ring whose radical is the fundamental ideal $\omega (G_0)$; the right socle of $k [G_0]$ is not zero since $k [G_0]$ is left perfect
and $G_0$ is finite \cite[p.~137]{ribenboim1969rings}, which contradicts the hypothesis made on $G$.

As a special case, we obtain the characterization of Artinian group rings
[I. G. Connell \cite{connell1963group}].

\section{Self-injective group rings}

\begin{theorem}\label{thm:4} Let $A$ be a ring and $G$ be a group. The following conditions are equivalent:

\begin{enumerate}
\item The ring $A [G]$ is left self-injective.
\item
\begin{enumerate}\item $A$ is left self-injective;
\item $G$ is a finite group.
\end{enumerate}
\end{enumerate}
\end{theorem}

(2) $\implies$ (1): This is a result of I. G. Connell \cite{connell1963group}.

(1) $\implies$ (2): Following \cite{connell1963group} we know that $A$ is left self-injective. Let $H$
be a finitely generated subgroup of $G$, and $\omega (H)$ be the right ideal of $A[G]$ generated
by the elements  $1 - h$, $h \in H$. According to \cite{ikdeda1954quasi} we know the left annihilator
of $\omega (H)$ is different from $(0)$, so $H$ is finite \cite{lambek1966lectures}, which proves
that $G$ is locally finite.

Suppose that $G$ is an infinite group; following \cite{hall1964property}, $G$ contains an infinite Abelian subgroup $G_1$. $A [G]$ which is a free $A [G_1]$ -module, is an
injective $A [G_1]$ -module \cite[p.~123]{cartan1960homological}, in particular $A [G_1]$ is left self-injective.
 If $H_1$ is an infinite subgroup of $G_1$, $A [G_1]$ is an injective $A [H_1]$ -module,
but as $A [H_1]$ is not a quasi-Frobenius ring
(See Theorem 3), this implies according to C. Faith \cite{faith1966rings}, that the index of $H_1$
in $G_1$ is finite. We deduce that the socle of $G_1$ is of finite length and $G_1$
is an Artinian Abelian group \cite{fuchs1968abelian}. It is easy to see the problem is reduced to
the case when $G_1$ is quasi-cyclic $p$-group. A contradiction results
from the following proposition [cf. also \cite{malliavin1971anneaux}].

\begin{proposition}\label{prop:5}(Pascaud).  Let $A$ be a ring and $G$ be the quasi-cyclic $p$-group
defined by generators $x_i$ and relations $ x_i = x ^p_{i + 1} $. Then $A [G]$
is not left self-injective. 
\end{proposition}

$A [G]$ is a free left $A$-module and we give $ B = Hom_A (A [G],
A [G]) $ a left $A [G]$ -module structure by  definining $ x \in A [G] $,
$ f \in  B$, $ (x\cdot f) (y) = f (yx)$ for $ y \in A [G]$.

$A [G]$ embeds into $B$ in the following way: to each $ x =\sum_ {g_i}a(g_i)  g_i$
we associate the endomorphism $ \bar x $: defined by $ \bar {x} (g_i) = a (g_i^{- 1})$.

We denote by $G_i $  the group generated by $x_i$ and we consider the elements $f$, $f_i$
of $B$ defined by:

$$
f(g)=\begin{cases} 
      1 & \text{ if } g=x^l_{2k}x_{2k+1} \text{ for some $k$, $l$} \\
      0 & \text{otherwise}
\end{cases}
$$

$$
f_i(g)=\begin{cases} 
      1 & \text{ if } g=x^l_{2k}x_{2k+1} \text{ for some $k$, $l$ with $k\leq i$} \\
      0 & \text{otherwise}
\end{cases}
$$

For all $i$, $f_i$ is an element of $A [G]$ and $f$ is an element of $B$ that
does not belong to $A [G]$.

\begin{lemma}\label{lem}(1) Let $a, b$ be two elements of $A [G_i]$,  $x$ an element of $G$,
$ x \notin  G_i$. The relation $a = bx $ implies $a = b = 0 $. 
(2) If $g$ is an element of $G$ not belonging to $ G_ {2i + 2} $, then
$ (1 - x_ {2i + 2})\cdot f (g) =  0$.
\end{lemma}

The proposition will result from the fact that $A [G] + A [G] f$ is an essential extension
of $A [G]$. Let $a, b$ two elements of $A [G_{2i}]$ with $ a + bf \neq  0$
If $ g \notin G_ {2i + 2} $, the support of $ bg $ does not meet $ G_ {2i + 2} $ and according to the Lemma
$ (1- x_ {2i + 2}) bf (g) = 0 $ and consequently
$$
y = (1-x_{2i + 2}) (a + bf) = (1-x_{2i + 2}) (a + bf_i)
$$
belongs to $A[G] $.  If $ y = 0 $, according to the lemma we have $ a + bf_i =  0$, from which it follows that $ a + bf = b (f-f_i) $. Let $  n_0$ be the smallest integer $ \geq  i + 1$ such that we have
$ b (f_n-f_i) \neq  0$; showing, as before, that
$$
(1 -x_{2n_0 + 2}) (a + bf) = (1-x_ {2n_0 +2}) b (f_ { n_0} -f_i)
$$
which is an element $\neq 0$ in $ A [G] $ according to property (1) of the Lemma.

\vskip 0.25in
\rightline{89 avenue du Recteur-Pineau}

\rightline{ 86-Poitiers, Vienne\hskip0.35in}
%%
% requires a BiBTeX file sample.bib

\bibliography{bibliography}{}

\providecommand{\bysame}{\leavevmode\hbox to3em{\hrulefill}\thinspace}
\providecommand{\MR}{\relax\ifhmode\unskip\space\fi MR }
\providecommand{\MRhref}[2]{%
  \href{http://www.ams.org/mathscinet-getitem?mr=#1}{#2}
}
\providecommand{\href}[2]{#2}
\begin{thebibliography}{10}

\bibitem{bass1960finitistic}
Hyman Bass, \emph{Finitistic dimension and a homological generalization of
  semi-primary rings}, Transactions of the American Mathematical Society
  \textbf{95} (1960), no.~3, 466--488.

\bibitem{bjork1969rings}
Jan-Erik Bj{\"o}rk, \emph{Rings satisfying a minimum condition on principal
  ideals.}, Journal f{\"u}r die reine und angewandte Mathematik \textbf{236}
  (1969), 112--119.

\bibitem{burgess1969semi}
Walter~D. Burgess, \emph{On semi-perfect group rings}, Canad. Math. Bull
  \textbf{12} (1969), no.~5, 645.

\bibitem{cartan1960homological}
Henri Cartan and Samuel Eilenberg, \emph{Homological algebra, princeton, 1956},
  Mathematical Reviews (MathSciNet): MR17: 1040e Zentralblatt MATH \textbf{75}
  (1960).

\bibitem{connell1963group}
Ian~G. Connell, \emph{On the group ring}, Canad. J. Math \textbf{15} (1963),
  no.~49, 650--685.

\bibitem{faith1966rings}
Carl Faith, \emph{Rings with ascending condition on annihilators}, Nagoya
  Mathematical Journal \textbf{27} (1966), no.~1, 179--191.

\bibitem{fuchs1968abelian}
L{\'a}szl{\'o} Fuchs, \emph{Abelian groups}, Budapest Hungarian Academy of
  Sciences, 1968.

\bibitem{gulliksen1970elementary}
Tor Gulliksen, Paulo Ribenboim, and TM~Viswanathan, \emph{An elementary note on
  group rings}, J. Reine Angew. Math \textbf{242} (1970), 148--162.

\bibitem{hall1964property}
Philip Hall and C.~R. Kulatilaka, \emph{A property of locally finite groups},
  Journal of the London Mathematical Society \textbf{1} (1964), no.~1,
  235--239.

\bibitem{herstein1971noncommutative}
Israel~Nathan Herstein, \emph{Noncommutative rings}, Mathematical Association
  of America, 1971.

\bibitem{ikdeda1954quasi}
Masatoshi Ikeda and Tadashi Nakayama, \emph{On some characteristic properties
  of quasi-frobenius rings and regular rings.}, Proceedings of the American
  Mathematical Society \textbf{5} (1954), 15--19.

\bibitem{lambek1966lectures}
Joachim Lambek, \emph{Lectures on rings and modules}, Blaisdell, 1966.

\bibitem{malliavin1971anneaux}
Marie-Paule Malliavin, \emph{Sur les anneaux de groupes fp-self-injectifs}, CR
  Acad. Sc. Paris Ser. A \textbf{273} (1971), 88--91.

\bibitem{faculte1970seminaire}
Guy Renault, \emph{S{\'e}minaire d'alg{\`e}bre non commutative}, Universit{\'e}
  de Paris, Facult{\'e} des sciences d'Orsay., 1970.

\bibitem{ribenboim1969rings}
Paulo Ribenboim, \emph{Rings and modules}, Interscience, 1969.

\bibitem{woods1971perfect}
Sheila Woods, \emph{On perfect group rings}, Proceedings of the American
  Mathematical Society \textbf{27} (1971), no.~1, 49--52.

\end{thebibliography}
\bibliographystyle{amsplain}
\newpage
\pagenumbering{gobble}
\includegraphics[width=\textwidth]{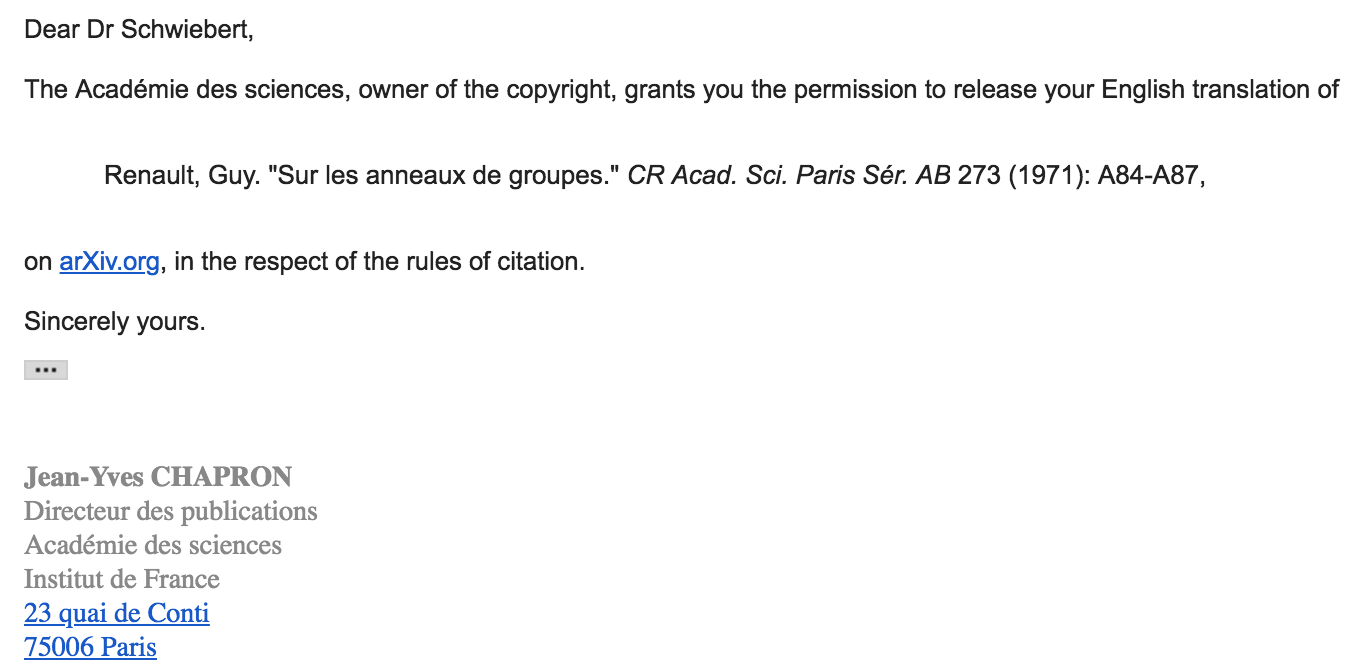}
\vskip 0.5in
The following changes were made to the original text owing to the high likelihood that they were typographical mistakes:
\begin{enumerate}
\item Page 1 third line of the introductory paragraph: ``Lambek'' was formerly ``Lambeck''.
\item Page 2 third line of intro to Section 2: ``Woods'' was formerly ``Wood''.
\item Page 2 line -2: The $k[G]$ at the beginning of the sentence was formerly $K[G]$.
\item Page 3 third line of Proposition 5: $ x =\sum_ {g_i}a(g_i)  g_i$ was formerly $ x \sum_ {g_i}a(g_i)  g_i$.
\item Page 3 first case in definition of $f$: the $x_{2k+1}$ was formerly $X_{2k+1}$.
\item Page 4 line 2: $a+bf\neq 0$ was formerly $a=bf\neq 0$

\end{enumerate}

\end{document}